\documentclass[12pt,reqno]{amsart}
\usepackage{amsmath,amssymb,amsfonts,amsthm}
\usepackage[mathscr]{eucal}
\usepackage{hyperref}

\author{Hovhannes Khudaverdian}
\author{Theodore Voronov}

\thanks{Based on a plenary talk given by Th.V.
at the XXII Bia{\l}owie{\.z}a meeting, June 2003.}
\address{Department of Mathematics, University of Manchester Institute of Science and Technology
(UMIST), United Kingdom}

\email{khudian@umist.ac.uk}
\email{theodore.voronov@umist.ac.uk}

\title[Odd Laplace operators and homotopy algebras]{Geometry of
differential operators, odd Laplacians, and homotopy algebras}

\newtheorem{thm}{Theorem}

\theoremstyle{definition}
\theoremstyle{remark}
\newtheorem{ex}{Example}[section]

\renewcommand{\leq}{\leqslant}

 \DeclareMathOperator{\ord}{ord}
 \DeclareMathOperator{\grad}{grad}
 \DeclareMathOperator{\Res}{Res}
 \DeclareMathOperator{\Ker}{Ker}
 \DeclareMathOperator{\im}{Im}

\newcommand{\Act}{{\mathscr{A}}}

\newcommand{\bro}{{\boldsymbol{\rho}}}
\newcommand{\brro}{{{\boldsymbol{\rho}}{\mathrm{'}}}}
\DeclareMathOperator{\Vol}{Vol} 
 
\DeclareMathOperator{\sub}{sub}

\newcommand{\lie}[1]{{\mathcal L}_{{#1}}}
\newcommand{\liex}[1]{{\lie{{#1}}}}
\renewcommand{\div}{\mathop{\mathrm{div}}}

\newcommand{\volx}{{Dx}}
\newcommand{\der}[2]{{\frac{\partial {#1}}{\partial {#2}}}}

\newcommand{\RR}{\mathbb R}
\newcommand{\Z}{{\mathbb Z_{2}}}

\newcommand{\p}{\partial}
\newcommand{\fun}{C^{\infty}}
\newcommand{\V}{{\mathfrak{V}}}

\def\a{\alpha}
\def\e{\varepsilon}
\def\s{\sigma}

\def\D{\Delta}

\newcommand{\g}{{\gamma}}

\renewcommand{\S}{{{S}}}
\newcommand{\lt}{\theta} 
\renewcommand{\r}{{\rho}}

\newcommand{\ft}{{\tilde f}}

\newcommand{\at}{{\tilde a}}
\newcommand{\bt}{{\tilde b}}

\newcommand{\ps}{{\boldsymbol{\psi}}}
\newcommand{\ch}{{\boldsymbol{\chi}}}

\newcommand{\DD}[1]{\Delta_{#1}}
\DeclareMathOperator{\Dr}{\DD{{\bro}}}
\DeclareMathOperator{\divr}{{\mathrm{div}_{\bro}}}
\newcommand{\Drr}{\DD{\brro}\,}
\DeclareMathOperator{\DDr}{\DD{{\bro}}^2}
\newcommand{\DDrr}{\DD{\brro}^2\,}

\begin{document}
\begin{abstract}
\noindent We give a complete description of differential operators
generating a given bracket.  In particular we consider the case of
Jacobi-type identities for odd operators and brackets. This is
related with homotopy algebras using the derived bracket
construction.
\end{abstract}

\maketitle 
\section{Introduction}

In this paper we give a survey of results of our
works~\cite{tv:laplace1}, \cite{tv:laplace2}, \cite{tv:higherder}
(see also \cite{tv:laplace2bis}). Their original motivation was to
give a clear geometric picture of the relation between an odd bracket
and its ``generating operator'' in the Batalin--Vilkovisky formalism.

In the pioneer paper~\cite{bv:perv}, for the needs of the Lagrangian
quantization of gauge theories, Batalin and Vilkovisky constructed a
remarkable second-order operator on the ``phase space''  of fields
and antifields.  One of us (H.K.) suggested a mathematical framework
in which the Batalin--Vilkovisky operator is interpreted as an
invariant operator acting on functions on an odd Poisson manifold
equipped with a volume form~\cite{hov:deltabest}. It is nothing but
an ``odd Laplace operator'' associated with an odd Poisson structure
in the same way as the usual Laplacian on a Riemannian manifold is
associated with a Riemannian metric, the main difference being in the
fact that on a Riemannian manifold there is a natural volume form,
while for an odd Poisson case, even if the bracket is non-degenerate
(i.e., for odd symplectic manifolds), no natural volume form exists
and a volume form should be introduced as an extra piece of data. (In
a similar situation in the even Poisson geometry a vector field
rather than a second-order operator arises~\cite{weinstein:modular}.)
As it has turned out, the bracket of functions can be recovered from
the odd Laplacian, as the failure of the Leibniz property, so the odd
Laplacian `generates' the odd Poisson structure.  Other constructions
of a generating operator, based on a connection, have been
proposed~\cite{yvette:divergence}. On the other hand, physical
motivations require  an operator acting on half-densities
(semidensities) rather than on functions. A canonical operator on
half-densities on an odd symplectic manifold was discovered
in~\cite{hov:max, hov:semi}. It does not require any extra structure
such as a connection or volume form. As it has turned out, on general
odd Poisson manifolds the situation is more
complicated~\cite{tv:laplace1}. There is no longer a unique operator
on half-densities, though the construction of a Laplace operator
using a volume form  is still `more canonical' for half-densities
than for functions or densities of any other weight, since for
half-densities this operator, as was shown in~\cite{tv:laplace1},
depends only on the orbit of a volume form under the action of a
certain groupoid. This groupoid (called the ``master groupoid''
in~\cite{tv:laplace1}) consists of transformations of volume forms
$\bro \mapsto \brro=e^f\bro$ satisfying the ``master equations'' $\Dr
e^{f/2}=0$, thusly revealing their composition property. (One should
note preliminary results in this direction in~\cite{ass:symmetry},
which in a hindsight can be interpreted as pointing at
half-densities.)

The results of~\cite{tv:laplace1} quite unexpectedly showed a
similarity between odd Poisson geometry and the usual Riemannian
geometry. Certain formal properties of the Laplacian on
half-densities or densities of other weight, happen to be the same
regardless of parity of the `bracket' (if one views a Riemannian
metric as a bracket). Physically, this points at a formal similarity
between the BV master equation on half-densities and the Schroedinger
equation, the classical master equation being analoguous to the
eikonal or Hamilton--Jacobi equation~\cite{tv:laplace1}. In fact, for
the four theories (odd/even Poisson, odd/even Riemannian structure)
defined by a rank two tensor $T^{ab}$ on $M$, analogies come in
pairs, see~\cite{tv:laplace1}.

Such a viewpoint helps, in a way, to `trivialize' the original
problem of the relation between a bracket and a generating operator,
by seeing it as the problem of describing all second-order operators
with a given principal symbol. (No confusion should be with the
quantization problem: the question is about one individual symbol,
not about a construction for all symbols.) In such form, a solution
is immediate. One only has to formulate it geometrically. For
second-order operators on functions, a piece of information necessary
to recover an operator from its principal symbol is the so called
subprincipal symbol. A good way of looking on it is to view it as a
sort of  connection, more precisely, an ``upper connection'' or
``contravariant derivative'' on volume forms, i.e., an operator
mapping them to vector fields. In a `bracket setup', it can  also be
viewed as an extension of the `bracket' of functions to a `bracket'
between functions and volume forms, or, using the Leibniz rule,
between functions and densities of arbitrary weight. This gives a
complete description of generating operators acting on functions.

While complete, this solution is not entirely satisfying
aesthetically, since it lacks symmetry: to describe operators just on
functions one has to consider brackets involving densities. Also,
operators, say, on half-densities must have their part in the
picture.

As we showed in~\cite{tv:laplace2}, passing from the algebra of
functions to the \textit{algebra of densities} $\V(M)$ on a
(super)manifold $M$ solves all problems. The reason for this is a
natural invariant scalar product in the algebra $\V(M)$.  It is
possible to establish a natural one-to-one correspondence between
`brackets' in this algebra and second-order operators. From the
viewpoint of $M$, a second-order operator in $\V(M)$ is a quadratic
pencil of operators $\D_w$ acting on densities of weight $w\in\RR$. A
bracket in $\V(M)$ is specified by a bracket of functions, a
corresponding to it upper connection on volume forms, and one extra
piece of data similar to the familiar to physicists Brans--Dicke
field of the Kaluza--Klein type models. These results contain all
previous formulae, e.g. for functions and half-densities, and explain
them.

As we said, up to a certain point there is no difference whatsoever
between the case of even operators and even brackets (which can be
considered on a usual manifold) and that of odd operators and odd
brackets, necessarily on a supermanifold. Parity of a bracket or of
its generating operator becomes essential when one wants to obtain
something like a Lie algebra. One can show that there is no way of
having a Jacobi-type identity for a symmetric even bracket
(`brackets' coming from operators are always symmetric, possibly in a
graded sense, because they are just polarizations of  quadratic
forms, the principal symbols). Therefore to make a progress in this
direction one should focus on the odd case, thus bringing us back to
`Batalin--Vilkovisky operators'. Our technique allows to give an
exhaustive answer to the questions about the possible ``Jacobi
identities'' involving an odd operator together with the
corresponding bracket.

From the algebraic viewpoint, identities that arise here are a
particular case of more general identities. This prompts to
investigate the situation further. The bracket generated by a
second-order operator is an example of  `derived brackets'
(see~\cite{yvette:derived}; for a recent survey
see~\cite{yvette:derived2}). Taken with the generating operator
itself, it should be viewed as a part of a sequence of $n$-ary
`brackets', which can be defined for any operator of order $N$, so
that $n=0,\ldots,N$. (This works for even and odd operators, and this
sequence of higher brackets has the property that the $k+1$ bracket
is equal to  the discrepancy of the Leibniz identity for the $k$-th
bracket, see~\cite{koszul:crochet85}.) If one asks for suitable
Jacobi identities, it is natural do it in an abstract setup. Such a
setup, called ``higher derived brackets'', has been suggested
in~\cite{tv:higherder}. We have shown there how simple data such as a
Lie superalgebra with a projector on an Abelian subalgebra and an odd
element $\D$ can produce a strongly homotopy Lie algebra; in fact, we
show how the `higher Jacobi identities' are controlled by $\D^2$. The
situation with differential operators and the corresponding brackets
is just one particular instance of that.

The algebraic constructions of~\cite{tv:higherder} make the first
step of generalizing the results of~\cite{tv:laplace2} to higher
order operators. This should be a subject of further studies.

This paper is organized as follows. In Section~\ref{sec:problem} we
formulate the problem and give examples. In particular we review the
properties of the operator $\Dr$.  In Section~\ref{sec:answer} we
introduce the algebra of densities $\V(M)$ and state the main theorem
about the one-to-one correspondence between operators and brackets.
In Section~\ref{sec:oddcase} we consider the case of odd operators
and brackets. In Section~\ref{sec:higher} we give the main results
concerning higher derived brackets. Throughout the paper we use the
supermathematical conventions about commutators, derivations, etc.;
tilde is used to denote parity.

We want to thank the organizers of the XXII Workshop on Geometric
Methods in Physics   for the invitation and for a fantastic
atmosphere at Bia{\l}owie{\.z}a during the meeting. We particularly
want to thank Anatole Anatolievich Odzijewicz and Mikhail
Aleksandrovich Shubin.

\section{Main problem: operators and brackets} \label{sec:problem}

\subsection{Setup}

Let $M$ be a supermanifold and $\D$ an arbitrary second-order
differential operator acting on $\fun(M)$. In local coordinates:
\begin{equation} \label{eq:operator}
\D=\frac{1}{2}\,S^{ab}\,\p_b\p_a+ T^a\,\p_a +R.
\end{equation}
The principal symbol of  $\D$ is the symmetric tensor field $S^{ab}$,
or the quadratic function $S=\frac{1}{2}\,S^{ab} p_bp_a$ on $T^*M$.
The principal symbol can be alternatively understood as a symmetric
bilinear operation on functions:
\begin{equation*}
\{f,g\}:=\D(fg)-(\D f)\,g-(-1)^{\e\ft}f\,(\D g) +\D(1)\,fg,
\end{equation*}
where $\e=\tilde\D$ is the parity of the operator $\D$. In
coordinates
\begin{equation*}
\{f,g\}=S^{ab}\p_bf\,\p_a g (-1)^{\at\ft}.
\end{equation*}
That $\{f,g\}$ is a bi-derivation can be formally deduced from the
fact that $\ord \D\leq 2$. In the following by a \textit{bracket} in
a commutative algebra we mean an arbitrary symmetric bi-derivation.
We say that $\D$ is a \textit{generating operator} for the bracket
$\{f,g\}$. Notice also that $\{f,g\}=[[\D,f],g]$, where $[\ ,\ ]$
denotes the commutator of operators.

\textit{Problem}:  construct a generating operator for a given
bracket, i.e., reconstruct $\D$ for a given symmetric tensor field
$S^{ab}$; describe all generating operators $\D$.

\subsection{Examples}

\begin{ex} Let
$(S^{ab})=(g^{ab})=(g_{ab})^{-1}$, where $(g_{ab})$ is a Riemannian
metric.  The Laplace--Beltrami operator
\begin{equation*}
\D f=\div \grad
f=\frac{1}{\sqrt{g}}\,\der{}{x^a}\left(\sqrt{g}g^{ab}\der{f}{x^b}\right),
\end{equation*}
where $g=\det(g_{ab})$, is  a generating operator for the ``bracket''
$\{f,g\}=\nabla f\cdot\nabla g=g^{ab}\p_a f\p_b g$ (up to a factor of
$\frac{1}{2}$). An arbitrary generating operator  for this bracket
has the appearance $\frac{1}{2}\D + X +R$, where $X$ is a vector
field, $R$ is a scalar function.
\end{ex}

For a degenerate $g^{ab}$, one has to replace $\sqrt{g}$ by some
volume density $\r$. In the following example this is always the
case, even for non-degenerate matrices.

\begin{ex} Suppose now $(S^{ab})$ is odd and specifies an odd Poisson
bracket. Take an arbitrary volume form $\bro=\r Dx$. Then the
operator introduced in~\cite{hov:deltabest},
\begin{equation}\label{eq:bvhov}
\Dr f=\divr X_f=\frac{1}{\r}\,\der{}{x^a}\left(\r
S^{ab}\der{f}{x^b}\right),
\end{equation}
is, up to $\frac{1}{2}$, a generating operator for the bracket. (We
denote by $X_f$  the Hamiltonian vector field corresponding to $f$.)
$\Dr$ mimics the Laplace--Beltrami operator.
\end{ex}

Unlike the classical Laplace--Beltrami above, the
operator~\eqref{eq:bvhov} requires a choice of $\bro$. We shall study
$\Dr$ in more detail below, as well as give a description of all
generating operators for a given odd Poisson bracket.

\begin{ex} Rather, a counterexample. If we try to apply the same
construction to an \textit{even} Poisson bracket, specified by a
Poisson tensor (bivector) $(P^{ab})$, which is antisymmetric, then
the second-order terms cancel, and we get a first-order operator
\begin{equation}\label{eq:modvect}
\Dr f=\divr X_f=\frac{1}{\r}\,\der{}{x^a}\left(\r
P^{ab}\right)\der{f}{x^b},
\end{equation}
known as a \textit{modular vector field} of an even Poisson bracket
(see~\cite{weinstein:modular}). It does not generate the bracket.
\end{ex}

\subsection{Properties of $\Dr$}

The operator $\Dr$ acting on functions on an odd Poisson manifold and
depending on a choice of volume form was introduced
in~\cite{hov:deltabest}.  (In~\cite{hov:deltabest}, everything was
formulated for odd symplectic manifolds, but the results hold for the
general Poisson case. See also~\cite{yvette:divergence}.) From now on
we redefine $\Dr$ by inserting the factor $\frac{1}{2}$ into
formula~\eqref{eq:bvhov}. The properties of $\Dr$ are as follows.
First,
\begin{equation*}
    \Dr (fg)= (\Dr f)g+(-1)^{\ft}f(\Dr g)+(-1)^{\ft+1} \{f,g\}\,.
\end{equation*}
If we change $\bro$, $\bro\mapsto \brro=e^\s \bro$, then
    $\Drr = \Dr +\frac{1}{2}\,X_{\s}$.
(Similar properties hold for an even  $S^{ab}=g^{ab}$.) The following
two properties are peculiar for an odd bracket. Holds
\begin{equation*}
    \Dr \{f,g\}= \{\Dr f,g\}+(-1)^{\ft+1}\{f,\Dr g\}.
\end{equation*}
The operator $\DDr$ is a Poisson vector field, and
    $\DDrr =\DDr- X_{H(\bro',\bro)}$,
where
    $H(\bro',\bro)= e^{-{{\s}}/{2}}\Dr(e^{{{\s}}/{2}})$.
See~\cite{tv:laplace1}. It follows that $\DDr$ gives a well-defined
cohomology class, which we call the \textit{modular class} for an odd
Poisson bracket.

\subsection{$\Dr$ on half-densities and the master groupoid}

The action of $\Dr$ extends to densities of arbitrary weight by
setting $\Dr\ps:=\bro^w\Dr(\bro^{-w}\ps)$ on $w$-densities. As it
turns out, the case of half-densities ($w=\frac{1}{2}$) is
distinguished. We have
\begin{equation*}
    [\Dr,f]= \liex{f}+\frac{1}{2}\,(1-2w)\,\Dr f
\end{equation*}
and
\begin{equation*}
    \D_{\bro'} = \Dr +\frac{1}{2}\,(1-2w)\liex{\s}-
    w(1-w)\,H(\bro',\bro)\,,
\end{equation*}
where $H(\bro',\bro)=e^{-{{\s}}/{2}}\Dr (e^{{{\s}}/{2}})$ as above
and $\liex{f}$ denotes the Lie derivative along a Hamiltonian vector
field $X_{f}$. Clearly $w=0,1,\frac{1}{2}$ are singular values. The
case $w=\frac{1}{2}$ is distinguished by a particularly simple
transformation law. On half-densities $\D_{\bro'} = \Dr
+\frac{1}{4}\,H(\bro',\bro)$. It follows that the solutions of the
`master equations' $\Dr e^{\s/2}=0$, for various $\bro$, can be
composed, making a groupoid, which we call the \textit{master
groupoid} of an odd Poisson manifold. On an orbit of the master
groupoid, the operator $\Dr$ acting on half-densities does not depend
on a volume form $\bro$.

In the odd symplectic case, one can show that the coordinate volume
forms corresponding to all Darboux coordinate systems belong to the
same orbit, giving a distinguished orbit \cite{hov:semi},
\cite{tv:laplace1}. We call this the `Batalin--Vilkovisky Lemma'
(compare~\cite{bv:closure}). It is a proper replacement of the
Liouville theorem, which is no longer valid in the odd case.
Therefore on odd symplectic manifolds, though there is no natural
volume form (and no such a form invariant under all canonical
transformations can exist), there exists a canonical $\D$-operator on
half-densities independent of a choice of volume form.

For operators on functions, the square $\DDr$, which is a Poisson
vector field in general, also does not change on an orbit of the
master groupoid. In particular, if for some $\bro$ it happens that
$\DDr=0$, then $\DDrr=0$ for all $\brro=e^\s\bro$ such that $\Dr
e^{\s/2}=0$. In the symplectic case, such is the ``Darboux coordinate
orbit''. In the general odd Poisson case, the existence of such an
orbit with $\DDr=0$ is an open question.

\section{Operators in the algebra of densities} \label{sec:answer}

\subsection{Subprincipal symbol as upper connection}

Come back to a general operator~\eqref{eq:operator} acting
\underline{on functions}. We can set $R=\D 1$ to zero without loss of
generality. Thus a reconstruction of $\D$ amounts to recovering the
coefficients $T^a$. Recall the notion of {H\"ormander}'s subprincipal
symbol; in our case $\sub\D=(\p_bS^{ba}(-1)^{\bt(\e+1)}-2T^a)p_a$.
Unlike the principal symbol, $\sub \D$ is coordinate-dependent.
Precisely, $\g^a=\p_bS^{ba}(-1)^{\bt(\e+1)}-2T^a$ has the
transformation law
\begin{equation*}
\g^{a'}=\left(\g^a+S^{ab}\,\p_b\ln J\right)\der{x^{a'}}{x^a}\,,
\end{equation*}
where $J=\frac{Dx'}{Dx}$ is the Jacobian.
Thus $\sub \D$ can be interpreted as an ``upper connection'' in the
bundle $\Vol M$, i.e., specifying a \textit{contravariant derivative}
$\nabla^a\rho=(S^{ab}\p_b +\g^a)\rho$ on volume forms. The
coordinate-dependent Hamiltonian $\g=\sub\D=\g^a p_a$ plays the role
of a local connection form. If the matrix $S^{ab}$ is invertible,
then we can lower the index $a$ to get a usual connection. (Notice
that  $\D$ acts on scalar functions, and no extra structure is
assumed on our manifold a priori. Geometry arises just from the
operator $\D$.) Thus, a second-order operator $\D$ on functions
(normalized by $\D 1=0$) is equivalent to a set of data: a bracket on
functions and an associated upper connection in $\Vol M$.

Alternatively,  such an upper connection in the bundle $\Vol M$ can
be viewed  as   an extension of the bracket of functions $\{f,g\}$ to
a `long' bracket $\{f,\ps\}$ where the second argument is a volume
form: $\{f,\ps\}= \left(S^{ab}\p_b f\,\p_a \psi(-1)^{\at\ft}
+\g^a\p_a f\, \psi\right)\volx$. Thus, $\D$ on functions is
equivalent to a bracket on functions equipped with an extension of it
to volume forms (for one argument). Now we want to make the situation
more symmetric.

\subsection{Interlude: the algebra of densities}

For a (super)manifold $M$ the \textit{algebra of densities} $\V(M)$
consists of formal linear combinations of densities of arbitrary
weights $w\in \RR$. It contains the algebra of functions. The
multiplication is the usual tensor product. We can specify elements
$\ps\in\V(M)$ by generating functions $\ps(x,t)$, which are defined
on the total space of a one-dimensional bundle $\hat M$ over $M$. The
algebra $\V(M)$ possesses a unit $1$ and a natural invariant scalar
product given by the formula: $\langle \ps,\ch\rangle=\int_M \Res
(t^{-2}\ps(x,t)\ch(x,t))\,Dx$. It can be viewed as a natural
generalized volume form on the supermanifold $\hat M$. Hence there is
a canonical divergence operator for derivations of the $\RR$-graded
algebra $\V(M)$.  Using it, one can classify derivations by
decomposing them into the divergence-free part and the ``scalar''
part (see~\cite{tv:laplace2}).

\subsection{Main theorem}

Consider brackets in the algebra $\V(M)$, i.e.,  symmetric
bi-derivations. A bracket of weight zero in $\V(M)$ is given by a
symmetric tensor $(\hat S^{\hat a\hat b})=
    \bigl(\begin{smallmatrix}
       S^{ab} & t  \g^a   \\
      t  \g^a & t^2\lt   \\
    \end{smallmatrix}
    \bigr)$ on $\hat M$. From the viewpoint of $M$, the blocks of
this matrix have the following meaning: $S^{ab}$ specifies a bracket
of functions, $\g^a$ gives an associated upper connection in $\Vol M$
as above, and a new bit of data $\lt$ allows to consider brackets
between volume forms, so that $\{Dx,Dx\}=\lt (Dx)^2$. It is a
third-order geometric object (the transformation law involves the
third derivatives of coordinates, see~\cite{tv:laplace2}, and depends
on $S^{ab}$ and $\g^a$), in the same way as $\g^a$ is a second-order
object with the transformation law depending on $S^{ab}$. So we have
something like a flag. The component $\lt$ is completely analogous to
the Brans--Dicke field $g^{55}$ in Kaluza--Klein type models.

On the other hand, let us consider differential operators in the
algebra $\V(M)$. This is stronger than simply take operators acting
on densities of various weights independently. In particular, a
second-order operator of weight zero in the algebra $\V(M)$  is a
quadratic pencil of operators on $w$-densities of the form
$\D_w=\D_0+wA+w^2B$ where $\D_0$ is a well-defined second-order
operator on functions,  $A$ and $B$ have orders one and zero,
respectively (they do not have invariant meaning separately from
$\D_0$). Since the algebra $\V(M)$ has a natural invariant scalar
product, it makes sense to consider self-adjoint operators. A pencil
$\D_w$ corresponds to a self-adjoint operator in $\V(M)$ if
$(\D_w)^*=\D_{1-w}$.

\begin{thm}
\label{thm:pencil} There is a one-to-one correspondence between
operators and brackets in $\V(M)$. Every second-order operator
generates a bracket and, conversely, for a given bracket a generating
operator always exists and can be uniquely specified by the
conditions of normalization $\D 1=0$ and self-adjointness. An
operator pencil $\D_w$ canonically corresponding to a bracket in
$\V(M)$ is given by the formula
\begin{equation*}
\D_w=\frac{1}{2}\left(\S^{ab}\p_b\p_a+
    \left(\p_b\S^{ba}(-1)^{\bt(\e+1)}+(2 w -1)\g^a\right)\p_a
    +\right. \\
    \hfill \left.
     w\,\p_a\g^a(-1)^{\at(\e+1)}  +
         w(w -1)\,\lt \right).
\end{equation*}
\end{thm}

The operator pencil defined in Theorem~\ref{thm:pencil} will be
shortly called the \textit{canonical pencil} (for a given bracket in
$\V(M)$).

\begin{ex}
The pencil considered in Section~\ref{sec:problem} and defined using
a volume form $\bro$, has $\g^a=S^{ab}\g_b$ and $\lt=\g^a\g_a$, where
$\g_a= \p_a\log \r$. Let us call it the \textit{Laplace--Beltrami
pencil}.
\end{ex}

The algebra $\V(M)$ allows to link  results for functions (when there
is a unit, but no scalar product) and for half-densities (when there
is a scalar product, but no unit). The canonical pencil corresponding
to a bracket in $\V(M)$ is nothing but the operator $\frac{1}{2}\div
\grad$ in $\V(M)$, where $\grad$ is given by a bracket and $\div$ is
the canonical divergence (see above). Using this description, one can
give the transformation law of the canonical pencil $\D_w$ under a
change of $\g^a$ and $\lt$ (see~\cite{tv:laplace2}), generalizing the
formula for $\Dr$. Taking the Laplace--Beltrami pencil as a
convenient `origin', one can get from there a useful parametrization
of all canonical pencils with a given $S^{ab}$,
see~\cite{tv:laplace2}. An interesting question is about the
specialization map $\D_w\mapsto \D_{w_0}$ from pencils to operators
on $w_0$-densities with a particular $w_0$. The values $w_0=0, 1,
\frac{1}{2}$ are singular. For other $w_0$, one can find a unique
canonical pencil such that $\D_{w_0}$ coincides with a given
second-order operator on $w_0$-densities (no restrictions). If
$w_0=\frac{1}{2}$, the image of the specialization map consists of
all self-adjoint operators on half-densities and the kernel consists
of pencils $(2w-1)\lie{X}$, where $X$ are vector fields. If $w_0=0$,
the image of the specialization map consists of all operators
vanishing on constants and the kernel is the subspace $\{w(w-1)f\}$.

\section{Odd case: Jacobi identities} \label{sec:oddcase}

\subsection{Algebraic statements}

Suppose $\D$ is an odd second-order differential operator in some
algebra $A$, generating an odd bracket in $A$, which we shall denote
$\{\ ,\ \}$. Set for simplicity $\D 1=0$ (we assume that there is a
unit). Notice that automatically $\ord \D^2\leq 3$, because
$\D^2=\frac{1}{2}[\D,\D]$.  It is not difficult to check the
following assertions:

$\ord\D^2\leq 2$ is equivalent to the Jacobi condition $\sum \pm
\{\{f,g\},h\}=0$; we shall refer to it as to
\underline{$\text{\textit{Jacobi}}_{\text{\,{3}}}$};

$\ord\D^2\leq 1$ is equivalent to the
\underline{$\text{\textit{Jacobi}}_{\text{\,{3}}}$} plus extra two
conditions, which are equivalent: $\D^2$ is a derivation of the
associative product, $\D$ is a derivation of the bracket; we shall
refer to the latter  as to
\underline{$\text{\textit{Jacobi}}_{\text{\,{2}}}$};

$\ord\D^2\leq 0$, finally,  is equivalent to all above plus $\D^2=0$,
to which we shall refer  as to
\underline{$\text{\textit{Jacobi}}_{\text{\,{1}}}$}.

Hence, in the notation
\underline{$\text{\textit{Jacobi}}_{\text{$\,{n}$}}$}, the number $n$
stands for the number of arguments.

\subsection{Geometric meaning}

We can apply this to operators and brackets in $\fun(M)$ or $\V(M)$.
In the latter case we know that there is a canonical operator
generating a given bracket (the canonical pencil). Keeping the above
notation, we have the following (parentheses stand for the canonical
Poisson bracket on $T^*M$):

\begin{thm}
For operators on functions:
\begin{equation*}
\ord \D^2\leq 2 \quad \Leftrightarrow\quad
\underline{\text{\textit{Jacobi}}_{\text{\,{3}}}} \quad
 \Leftrightarrow \quad  (S,S)=0,
\end{equation*}
hence $D=(S,\ )$ is a differential;
\begin{equation*}
\ord \D^2\leq 1 \quad \Leftrightarrow\quad
\underline{\text{\textit{Jacobi}}_{\text{\,{3}}}} \  + \
\underline{\text{\textit{Jacobi}}_{\text{\,{2}}}} \quad
 \Leftrightarrow \quad  (S,S)=0, \ (S,\g)=0,
\end{equation*}
i.e., $\g$ is flat, $D\g=0$.
\end{thm}

Notice that $D\g=(S,\g)$ plays the role of curvature for an upper
connection $\g$; it only makes sense with $D^2=0$, i.e., when $\{\ ,\
\}$ is a genuine odd Poisson bracket.

A further condition that $\D^2=0$ might be seen as a version of a
`Batalin--Vilkovisky equation' for a pair $S, \g$ of an odd Poisson
bracket and a flat upper connection, but we  prefer to relate `BV
equations' with changes of $\g$, like for $\bro$ in
Section~\ref{sec:problem}.

\begin{thm}
For a unique operator corresponding to a bracket in $\V(M)$: the
conditions
\begin{equation*}
\ord \D^2\leq 2 \quad \Leftrightarrow\quad
\text{\underline{$\text{\textit{Jacobi}}_{\text{\,3}}$} \, in\,\,
$\V(M)$}
\end{equation*}
automatically imply the conditions
\begin{equation*}
\ord \D^2\leq 1 \quad \Leftrightarrow\quad
\underline{\text{\textit{Jacobi}}_{\text{\,3}}} \,+\,
\underline{\text{\textit{Jacobi}}_{\text{\,2}}} \  \text{\,in\,\,
$\V(M)$}
\end{equation*}
and are equivalent to the equations
\begin{align*}
    (S,S)=0, \     (S,\g)=0    \\
    (S,\lt)+(\g,\g)=0  \\
    (\g,\lt)=0
\end{align*}
\end{thm}

We can get an interesting corollary from here. Suppose the bracket
generated on functions is non-degenerate (an odd symplectic
structure). Then the equations for $\g, \lt$ are solved uniquely (if
we ignore constants) giving $\g^a=S^{ab}\g_b$ where $\g_b=-\p_b
\log\Act=e^{\Act}\p_b e^{-\Act}$ for some volume form
$\bro=e^{-\Act}Dx$, where $\Act$ can be interpreted as ``action'',
and $\lt=\g^a\g_a$. Thus in this case the Jacobi conditions in the
algebra of densities bring us back to the operator $\Dr$.

\section{Homotopy algebras} \label{sec:higher}

\subsection{Brackets generated by operators of higher order}

The following construction was essentially given by
Koszul~\cite{koszul:crochet85}. Let $\D$ be an operator in an
arbitrary commutative $\Z$-graded algebra with a unit. Define
\begin{equation*}
    \begin{split}
    \{a\}&=[\D,a] 1,    \\
    \{a,b\}&=[[\D,a],b] 1, \\
    \{a,b,c\}&=[[[\D,a],b],c] 1, \\
    \text{etc.};
    \end{split}
\end{equation*}
the $n$-ary bracket is obtained by taking the $n+1$ consequent
commutators. One might add an $0$-ary bracket, which is simply the
element $\D 1$. $\D$ is a differential operator of order $\leq N$ if
all brackets with more than $N$ arguments vanish. In such case the
top nonzero bracket is a multi-derivation. It is the polarization of
the symbol of $\D$. For $N=2$ we return to the situation considered
above. One can check that all brackets obtained in this way are
symmetric and that the $(k+1)$-fold bracket appears as the
obstruction to the Leibniz rule for the $k$-fold bracket. This holds
for even and odd operators.

\subsection{Case of odd operators}

Suppose that $\D$ is odd. Hence all the brackets generated by $\D$
are also odd. Can we get, in addition to those linked `Leibniz
identities' holding automatically, some sort of Jacobi identities for
these brackets?

Let $a_1, \ldots, a_n$ be elements of $A$. We do not need
multiplication here, so for a moment one can think of $A$  just as of
a vector space. Define the \textit{$n$-th Jacobiator} on $a_1,
\ldots, a_n$ by the formula

\begin{equation*}
    J^n(a_1, \ldots, a_n)=
    \sum_{\parbox{1.2cm}{\small $\scriptstyle k+l=n$}} \sum_{\text{$(k,l)$-shuffles}}
    (-1)^{\a}
    \{\{a_{\s(1)},\ldots,a_{\s(k)}\},a_{\s({k+1})},\ldots,a_{\s({k+l})}\}.
\end{equation*}
Here $(-1)^{\a}$ is the sign prescribed by the sign rule for a
permutation of homogeneous elements $a_1,\ldots, a_n\in A$. This
definition assumes the existence  in $A$ of odd brackets of all
orders between $0$ and $n$. The \textit{$n$-th Jacobi identity} for a
sequence of odd brackets in $A$ is the  equality $J^n=0$ identically
for all arguments.

Up to sign conventions, this is essentially the definition used in
Stasheff's theory of strongly homotopy algebras: a vector space
endowed with a sequence of (super)symmetric odd  brackets such that
all Jacobi identities hold, for all $n=0,1,2,3,\ldots$,  is a
\textit{strongly homotopy Lie algebra} (an
\textit{$L_{\infty}$-algebra}). Very often it is assumed that the
distinguished element given by the $0$-ary bracket and called a
\textit{background} vanishes.  See, e.g., ~\cite{lada:stasheff},
\cite{lada:shla},  where  slightly different conventions are used.

Coming back to our situation when the brackets are generated by an
odd operator $\D$, we have the following remarkable
statement~\cite{tv:higherder}.

\begin{thm} \label{thm:higherbrack}
The following conditions are equivalent:

$(1)$ $\ord \D^2\leq r$;

$(2)$ all Jacobi identities with the number of arguments $>r$ are
satisfied.
\end{thm}

In particular, if $\D^2=0$, we get an $L_{\infty}$-algebra structure
in $A$  besides the associative multiplication (compare
also~\cite{bering:higher}). This is a particular example of
``homotopy Batalin--Vilkovisky algebras''.

Theorem~\ref{thm:higherbrack} is a corollary of the general
Theorem~\ref{thm:higherabstr} below,  valid in the following abstract
axiomatic setting suggested in~\cite{tv:higherder}. Consider a Lie
superalgebra $\mathfrak g$ endowed with a projector $P$ such that
$\im P$ and $\Ker P$ are subalgebras and $\im P$ is Abelian. Let
$D\in\mathfrak g$ be an arbitrary element. The \textit{$n$-th derived
bracket} of $D$ is a symmetric multilinear operation on the space
$V=\im P$ defined by the formula
\begin{equation*}
    \{a_1,\ldots,a_n\}_{\D}:=P[\ldots[[\D,a_1],a_2],\ldots,a_n],
\end{equation*}
where $a_i\in V$.

\begin{thm} \label{thm:higherabstr}
Consider an odd element $\D\in\mathfrak g$.  The $n$-th Jacobiator of
the derived brackets of $\D$ is equal to the $n$-th derived bracket
of the even element $\D^2$:
\begin{equation*}
    J^n_{\D}(a_1,\ldots,a_n)=\{a_1,\ldots,a_n\}_{\D^2}.
\end{equation*}
\end{thm}

In our example the Lie superalgebra $\mathfrak g$ consists of all
operators in a commutative associative algebra with a unit $A$, the
projector $P$ is the evaluation on the unit element, $\im P$ being
$A$. Theorem~\ref{thm:higherbrack} immediately follows.



\def\cprime{$'$} \def\cprime{$'$}

\end{document}